\documentclass[12pt]{amsart}
\usepackage{amssymb,amsmath,amsthm}
\usepackage{graphicx}
\usepackage{subcaption}
\usepackage{xcolor}
\def\F {{\mathsf F}}
\def\E {{\mathsf E}}
\def\G {{\mathsf G}}
\def\R {\mathbb{R}}
\def\H {\mathcal{H}}
\def\D {{\rm{dom}}}
\def\A {{\mathbb A}}
\def \au {\rm}
\def \ti {\it}
\def \jou {\rm}
\def \bk {\it}
\def \no#1#2#3 {{\bf #1} (#3), #2.}
\def \eds#1#2#3 {#1, #2, #3.}
\newtheorem{proposition}{Proposition}[section]
\newtheorem{theorem}[proposition]{Theorem}
\newtheorem{corollary}[proposition]{Corollary}
\newtheorem{lemma}[proposition]{Lemma}
\theoremstyle{definition}
\newtheorem{definition}[proposition]{Definition}
\newtheorem{remark}[proposition]{Remark}
\newtheorem{example}[proposition]{Example}
\numberwithin{equation}{section}
\title[Abstract damped wave equations]
{Abstract damped wave equations:\\
The optimal decay rate}

\author[F. Dell'Oro, L. Liverani, V. Pata]
{Filippo Dell'Oro, Lorenzo Liverani, Vittorino Pata}
\thanks{Work partially supported by the Italian MIUR-PRIN Grant 2020F3NCPX “Mathematics for industry 4.0 (Math4I4)”}
\address{Politecnico di Milano - Dipartimento di Matematica
\newline\indent
Via Bonardi 9, 20133 Milano, Italy}
\email{filippo.delloro@polimi.it {\rm (F. Dell'Oro)}}
\email{lorenzo.liverani@polimi.it {\rm (L. Liverani)}}
\email{vittorino.pata@polimi.it {\rm (V. Pata)}}

\subjclass[2010]{35B35, 35P05, 47D06}
\keywords{Abstract damped wave equations, semigroups, exponential stability, decay rate, resonance}
\begin{document}

\begin{abstract}
The exponential decay rate of the semigroup $S(t)=e^{t\A}$ generated by the
abstract damped wave equation
$$\ddot u + 2f(A) \dot u +A u=0
$$
is here addressed,
where $A$ is a strictly positive operator. The continuous function $f$, defined on the spectrum of $A$,
is subject to the constraints
$$\inf f(s)>0\qquad\text{and}\qquad \sup f(s)/s <\infty$$
which are known to be necessary and sufficient for exponential stability to occur.
We prove that the operator norm of the semigroup fulfills the estimate
$$\|S(t)\|\leq Ce^{\sigma_*t}$$
being $\sigma_*<0$ the supremum
of the real part of the spectrum of $\A$. This estimate always holds except in the resonant cases,
where the negative exponential $e^{\sigma_*t}$ turns out to be penalized by a factor $(1+t)$.
The decay rate is the best possible
allowed by the theory.
\end{abstract}

\maketitle
%%%%%%%%%%%%%%%%%%%%%%%%%%%%%%%%%%%%%%%%%%%%

%%%%%%%%%%%%%%%%%%%%%%%%%%%%%%%%%%%%%%%%%%%%
\section{Introduction}

\subsection{The equation}
Let $H$ be a (separable) real Hilbert space, endowed with the scalar product and norm
$\langle \cdot , \cdot \rangle,$ and $\|\cdot\|$, respectively, and let
$A :\D(A)\subset H \to H$
be a strictly positive selfadjoint linear operator.
We define the product Hilbert space
$$\H=\D(A^{\frac12})\times H,$$
endowed with the standard Euclidean product norm
$$\|(u,v)\|_\H^2=\|A^{\frac12}u\|^2+\|v\|^2.$$
Denoting by $\sigma(A)\subset(0,\infty)$ the spectrum of $A$, let
$$f:\sigma(A)\to(0,\infty)
$$
be a continuous function subject to the structural assumptions
\begin{equation}
\label{STRUCT}
\inf_{s\in\sigma(A)} f(s)>0\qquad\text{and}\qquad
\sup_{s\in\sigma(A)} \frac{f(s)}{s} <\infty.
\end{equation}
For $t>0$, we consider the abstract damped wave equation in the unknown $u=u(t)$
\begin{equation}
\label{WAVE}
\ddot u + 2f(A) \dot u +A u=0,
\end{equation}
the {\it dot} standing for derivative with respect to time,
subject to the initial conditions
$$u(0)=u_0\quad\text{and}\qquad \dot u(0)=v_0,$$
where the vector $(u_0,v_0)\in\H$ is arbitrarily assigned.
Here, $f(A)$ is the selfadjoint linear operator on $H$
constructed via the functional calculus of $A$,
namely,
$$
f(A) = \int_{\sigma(A)} f(s) \, d E_A (s),
$$
being $E_A$ the spectral measure of $A$ (see, e.g., \cite{RUD}).
Calling $\boldsymbol{u}=(u,v)$, and introducing the linear operator
$$
\A(u,v)=(v,- A u - f(A)v),$$
with domain
$$
\D(\A)=\big\{ \boldsymbol{u}\in\H:\,
\A\boldsymbol{u}\in\H\big\},
$$
we can rewrite~\eqref{WAVE}
as the first order ODE in $\H$
$$
\dot{\boldsymbol{u}}=\A \boldsymbol{u}.
$$
The linear operator $\A$ is known to be the infinitesimal generator
of a ${\mathcal C}_0$-semigroup
$$S(t)=e^{t\A}:\H\to\H$$
of linear contractions (see, e.g., \cite{DPWAVE}).
Hence, for any given $\boldsymbol{u}_0=(u_0,v_0)\in\H$, the unique mild solution
$\boldsymbol{u}(t)$ to the equation above, subject to the initial
condition $\boldsymbol{u}(0)=\boldsymbol{u}_0$,
is given by
$$\boldsymbol{u}(t)=(u(t),\dot u(t))=S(t)\boldsymbol{u}_0.
$$
We denote (twice) the energy of such a solution by
$$\E(t)=\|\boldsymbol{u}(t)\|^2_\H=\|A^{\frac12}u(t)\|^2+\|\dot u(t)\|^2.$$

\subsection{Exponential decay of the semigroup}
As far as the sole existence of $S(t)$ is concerned,
we could have taken
\emph{any} continuous function $f:\sigma(A)\to [0,\infty)$. However,
as shown for instance in~\cite{DPWAVE},
the exponential stability of the semigroup occurs if and only if~\eqref{STRUCT} is in place.
Recall that $S(t)$ is said to be \emph{exponentially stable} if
there exist $\omega>0$ and $C\geq 1$ such that
\begin{equation}
\label{expo}
\|S(t)\|_{L(\H)} \leq C e^{-\omega t},
\end{equation}
where the norm is taken in the Banach space $L(\H)$ of bounded linear operators on $\H$.
Once the exponential stability is known,
one would like to find the best possible $\omega>0$
(i.e., the largest) for which~\eqref{expo} holds. This number is the (exponential) \emph{decay rate}, defined
as\footnote{The number $-\omega_*$ is usually called the \emph{growth bound} of the semigroup.}
$$
\omega_*=\sup\big\{\omega > 0:\,\|S(t)\|_{L(\H)}\leq  C e^{-\omega t}\big\},
$$
for some $C=C(\omega)\geq1$.
However, for general (exponentially stable) semigroups, computing $\omega_*$ might not be an easy task.
Another relevant quantity, usually much easier to detect, is
the \emph{spectral bound} of $\A$
$$
\sigma_*=\sup_{\lambda\in\sigma(\A)}{\mathfrak{Re}}\, \lambda,
$$
$\sigma(\A)$ being the spectrum of (the complexification of) $\A$,
which is related to $\omega_*$ through
the (possibly strict) inequality $\omega_*\leq -\sigma_*$ (see \cite{ENG,PAZ}).

\begin{definition}
The semigroup $S(t)$ satisfies the \emph{spectrum determined growth} (SDG) condition
if $\omega_*=-\sigma_*$.
\end{definition}

In view of finding the decay rate, it is then
of paramount importance to establish whether or not the SDG condition is satisfied.
This is true, for instance, for eventually norm continuous (such as analytic or differentiable)
semigroups~\cite{ENG}. Still, even if $S(t)$ fulfills the SDG condition, this does not mean at all that~\eqref{expo} holds for
$\omega=\omega_*$, as the following simple example shows.

\begin{example}
\label{expendulum}
Consider the damped pendulum equation
$$\ddot u + 2a \dot u +u=0,$$
a particular instance of~\eqref{WAVE} for $\H=\R^2$, $A=1$
and $f\equiv a>0$.
The related semigroup, being norm continuous, always fulfills the SDG condition. In particular,
the decay rate is attained whenever $a\neq 1$. But, when $a=1$,
the norm of the solution corresponding to the initial datum $\boldsymbol{u}_0=(u_0,v_0)$ reads
$$\|S(t)\boldsymbol{u}_0\|_{\H}=\sqrt{u_0^2+v_0^2+2(u_0^2-v_0^2)t+2(u_0+v_0)^2t^2}\,e^{-\frac{t}2}.
$$
The decay rate, equal to $\frac12$ in this case, is certainly not attained, and the reason is that
a resonance phenomenon is encountered.
\end{example}

This motivates the next definition.

\begin{definition}
The semigroup $S(t)$ satisfies the \emph{strong spectrum determined growth} (SSDG) condition
if the decay rate $\omega_*=-\sigma_*$ is attained, that is, if~\eqref{expo} holds with
$\omega=\omega_*=-\sigma_*$.
\end{definition}

\subsection{Earlier contributions}
The equation~\eqref{WAVE} of this paper belongs to a general class of abstract differential
models introduced in~\cite{CR}, and ``exhibiting the empirically observed damping rates in elastic systems".
For this reason, $A$ is usually called the \emph{elastic operator},
while $f(A)$ is the \emph{damping operator}.
The prototypical case occurs when $f(s)$ controls
and is controlled by $s^\theta$ for some $\theta \in [0,1]$,
namely, when $f(A)$ is comparable with the powers $A^\theta$. Here,
the semigroup is know to be analytic for
$\theta \in[\frac12,1]$ and differentiable (actually, of Gevrey type) for $\theta \in (0,\frac12)$.
See, among others, the classical papers~\cite{CT1,CT2,CT3,HUA1,HUA2,HUA3}. Further works,
dealing with a damping operator not necessarily comparable
with $A^\theta$, include~\cite{BATK,GRIN1,JT,LTBOOK,LIU}.
Nevertheless, the literature concerning
the analysis of the decay rate $\omega_*$ of $S(t)$, with particular reference
to whether or not it is actually attained,
is exceptionally bare.
The full solution to the problem has been found
only for the weakly damped wave equation,
sometimes referred to as the \emph{telegrapher's equation},
corresponding to the choice $f(s)=a>0$ (see~\cite{SILVER,GOLD1}).
The much more general (and difficult) situation of a nonconstant $f$ has been
tackled in~\cite{GOLD2}. There, within a certain number
of assumptions, which prevent $f(s)$ to grow at infinity
faster than $s^\theta$ with $\theta<\frac12$, the authors obtain sharp exponential decay estimates
for trajectories originating from sufficiently regular initial data, but they cannot generally
exhibit the decay rate of the semigroup.
On the other hand, the necessary
and sufficient condition~\eqref{STRUCT} for the exponential stability of $S(t)$ clearly
allows $f(s)$ to have up to a linear growth
at infinity. And indeed, the special case $f(s)=s$ corresponds to the widely
studied strongly damped wave equation,
also known by the name of \emph{Kelvin-Voigt equation}, which appears in several areas of Mathematical Physics.

\subsection{Our result}
The aim of this paper is to provide a complete answer.
Within the sole
condition~\eqref{STRUCT}, we show that the semigroup $S(t)$
fulfills~\eqref{expo} with
$\omega=\omega_*=-\sigma_*$, except in some particular cases, called resonant,
where the term $e^{-\omega_*t}$ is penalized by a factor $(1+t)$.
This result is optimal. The key idea of the proof
is to decompose $\sigma(A)$
into the disjoint union of suitably chosen subsets $\sigma_\imath$, and then perform energy estimates
with \emph{ad hoc} multipliers
on each of the (orthogonal) subspaces $E_A(\sigma_\imath)H$,
making use of the functional calculus of $A$.

The paper is organized as follows. In the next Section~\ref{Sec2}, we introduce
two objects that will play a crucial role in our analysis.
In Section~\ref{Sec3}, we recall the properties of the spectrum of the operator $\A$.
The main result is stated and proved in Section~\ref{Sec4}, and an application is then discussed
in Section~\ref{Sec5}.
The final Sections~\ref{SecL0}-\ref{SecL3} are devoted to the proofs of four lemmas
encountered in the proof of the main theorem.

\subsection*{A word of warning}
Although we work with a real Hilbert space $H$, having in mind the concrete examples of
Mathematical Physics,
all the results of this paper hold \emph{verbatim} if $H$ is a complex Hilbert space.
The proofs are exactly the same, up to replacing in the calculations any occurrence of a scalar product
with its real part.
%%%%%%%%%%%%%%%%%%%%%%%%%%%%%%%%%%%%%%%%%%%%

%%%%%%%%%%%%%%%%%%%%%%%%%%%%%%%%%%%%%%%%%%%%
\section{The Function ${\phi}$ and the Number ${m_*}$}
\label{Sec2}

\noindent
A crucial object for our analysis is the continuous function $\phi:\sigma(A)\to(0,\infty)$
defined as
$$
\phi(s)=\begin{cases}
f(s) &\quad\text{if } f(s)\leq \sqrt{s}\,,\\
\noalign{\vskip1.5mm}
f(s) -\sqrt{f^2(s)-s} &\quad\text{if } f(s)> \sqrt{s}\,,
\end{cases}
$$
along with the number
$$
m_*=\inf_{s\in \sigma(A)} \phi(s)>0.
$$
The fact that $m_*>0$ merely follows by observing that
$$
f(s)=\begin{cases}
\phi(s) &\quad\text{if } f(s)\leq \sqrt{s}\,,\\
\noalign{\vskip.5mm}
\displaystyle\frac{s}{2\phi(s)}+\frac{\phi(s)}{2}&\quad\text{if } f(s)> \sqrt{s}\,.
\end{cases}
$$
Hence, $m_*=0$ would violate~\eqref{STRUCT}.
Then, in the region $f(s)> \sqrt{s}\,$, we have that
\begin{equation}
\label{ossbella}
\frac{f(s)}{s}\leq \frac1{2m_*}+\frac{m_*}{2s},
\end{equation}
this fact being equivalent to $\phi(s)\geq m_*$.
Finally, for every $m>0$ and every sequence $s_n\to\infty$, we have the implication
\begin{equation}
\label{liminfi}
\frac{f(s_n)}{s_n}\to\frac{1}{2m}
\quad\Rightarrow\quad \phi(s_n)\to m.
\end{equation}
Indeed, the first convergence implies in particular that
$$\frac{s_n}{f^2(s_n)}\to0\qquad\text{and}\qquad f(s_n)>\sqrt{s_n}\,,$$
the latter for all $n$ large. Accordingly,
$$\phi(s_n)=f(s_n)\Big(1-\sqrt{1-\frac{s_n}{f^2(s_n)}}\,\Big)\sim\frac{s_n}{2f(s_n)}\to m.
$$
%%%%%%%%%%%%%%%%%%%%%%%%%%%%%%%%%%%%%%%%%%%%

%%%%%%%%%%%%%%%%%%%%%%%%%%%%%%%%%%%%%%%%%%%%
\section{The Spectrum of the Infinitesimal Generator}
\label{Sec3}

\noindent
For every fixed $s\in\sigma(A)$, we introduce the pair of complex numbers
$$
\lambda_s^{\pm}=
\begin{cases}
\displaystyle -f(s) \pm i\sqrt{s-f^2(s)} &\quad\text{if } f(s)\leq \sqrt{s}\,,\\
\noalign{\vskip1mm}
\displaystyle-f(s)\pm  \sqrt{f^2(s)-s} &\quad\text{if } f(s)> \sqrt{s}\,,
\end{cases}
$$
which are nothing but the solutions to the second order equation
$$
\lambda^2+ 2f(s)\lambda +s=0.
$$
We define
$$\Sigma=\bigcup_{s\in\sigma(A)} \big\{\lambda_s^{\pm} \big\}.
$$
We also consider the (possibly empty) set
$$
\Lambda = \Big\{\lambda<0:\, \exists \, s_n\in\sigma(A)
\,\,\,\text{such that}\,\,\, s_n\to\infty \,\,\,\text{and}\,\,\, \lim_{n\to\infty} \frac{f(s_n)}{s_n}=-\frac{1}{2\lambda}\Big\}.
$$
Then, as shown in~\cite{DPWAVE}, the spectrum of (the complexification of) the operator $\A$
reads
$$
\sigma(\A) = \Sigma \cup \Lambda,
$$
where the union is not necessarily disjoint.
In light of~\eqref{liminfi}, it is apparent that $\Lambda$ belongs to the closure of $\Sigma$.
This yields the equality
\begin{equation}
\label{eqSpectralBound}
\sigma_*=\sup_{\lambda\in\sigma(\A)}{\mathfrak{Re}}\, \lambda
=\sup_{\lambda\in\Sigma}{\mathfrak{Re}}\, \lambda
=-\inf_{s\in \sigma(A)} \phi(s)=-m_*.
\end{equation}

\begin{remark}
Concerning the structure of the spectrum of $\A$, we have that
$$\sigma(\A)=\sigma_{\rm p}(\A)\cup\sigma_{\rm c}(\A),
$$
where
$$\sigma_{\rm p}(\A)=\bigcup_{s\in\sigma_{\rm p}(A)} \big\{\lambda_s^{\pm} \big\}.
$$
With standard notation, $\sigma_{\rm p}$ and $\sigma_{\rm c}$
are the point spectrum and the continuous spectrum, respectively.
This characterization, although not explicitly stated in~\cite{DPWAVE}, can be
inferred from the calculations of that paper, by arguing as in~\cite[Sec.\ 5]{VDRM}.
\end{remark}
%%%%%%%%%%%%%%%%%%%%%%%%%%%%%%%%%%%%%%%%%%%%

%%%%%%%%%%%%%%%%%%%%%%%%%%%%%%%%%%%%%%%%%%%%
\section{The Theorem}
\label{Sec4}

\subsection{Statement of the result}
We begin with the rigorous definition of resonance.

\begin{definition}
\label{RESDEF}
The semigroup $S(t)$ is said to be \emph{resonant} if
there exists $s_*\in\sigma(A)$ such that
$$m_*=\phi(s_*) \qquad\text{and}\qquad f(s_*)=\sqrt{s_*}\,.$$
\end{definition}

In other words, the resonance phenomenon appears whenever
the real part of the spectrum of $\A$ attains its supremum $-m_*$ and the equality
$$-m_*=\lambda^+_{s_*}=\lambda^-_{s_*}$$
holds for some $s_*\in\sigma(A)$. In fact, if such a point $s_*$ exists, it is easily seen to be unique.

\begin{theorem}
\label{MAIN}
There exists a constant $C\geq 1$ such that
\begin{itemize}
\smallskip
\item[{\bf I.}] $\|S(t)\|_{L(\H)} \leq C e^{-m_* t},\,\,\,$ if $S(t)$ not resonant; and
\medskip
\item[{\bf II.}] $\|S(t)\|_{L(\H)} \leq C (1+t) e^{-m_* t},\,\,\,$ if $S(t)$ resonant.
\end{itemize}
\smallskip
\end{theorem}

On account of~\eqref{eqSpectralBound}, the theorem produces a straightforward consequence.

\begin{corollary}
If $S(t)$ is not resonant, then it fulfills the {\rm SSDG} condition. Instead,
if $S(t)$ is  resonant, it fulfills the {\rm SDG} condition, but not the {\rm SSDG} one.
\end{corollary}

\begin{remark}
In view of Example~\ref{expendulum}, the result is optimal.
\end{remark}

It is sometimes possible to provide a precise estimate of the constant $C$.

\begin{proposition}
\label{propsub}
Assume that the semigroup is subdamped, that is,
$$
\sup_{s\in\sigma(A)}\frac{f(s)}{\sqrt{s}}=\ell<1.
$$
Then point {\bf I} of Theorem~\ref{MAIN} holds with
$$C=\sqrt{\frac{1+\ell}{1-\ell}}.$$
\end{proposition}

Also in this case, the result is optimal, as the following example shows.

\begin{example}
Consider the subdamped pendulum equation
$$\ddot u + 2 a\dot u +u=0,\quad a\in(0,1),$$
which meets the hypotheses of the proposition above with $\ell=a$.
Taking the unit vector $\boldsymbol{u}_0=(\frac1{\sqrt{2}}\,,\frac1{\sqrt{2}})$,
the norm of the solution at time $t$ reads
$$\|S(t)\boldsymbol{u}_0\|_{\H}
=\sqrt{\frac{1-a\cos\big(2\sqrt{1-a^2}\, t\big)}{1-a}}\,e^{-\frac{a}2t}.
$$
Choosing
$$t=t_k=\frac{(2k+1)\pi}{2\sqrt{1-a^2}\,},$$
for any $k\in{\mathbb N}$, we get exactly
$$\|S(t_k)\|_{L(\H)}\geq\|S(t_k)\boldsymbol{u}_0\|_{\H}
=\sqrt{\frac{1+a}{1-a}}\,e^{-\frac{a}2{t_k}}.$$
\end{example}

\subsection{Some preparatory lemmas}
For any given $K\geq 2$ and $\varepsilon\in(0,1)$, we decompose the spectrum of $A$ into the disjoint union
$$\sigma(A)=\sigma_0\cup\sigma_1\cup \sigma_2\cup \sigma_3,
$$
where the four regions $\sigma_\imath$ (some of which possibly empty) are defined as
\begin{align*}
\sigma_0&=\Big\{s\in\sigma(A):\,\frac{f(s)}{\sqrt{s}}>K\Big\},\\
\sigma_1&=\Big\{s\in\sigma(A):\,\frac{f(s)}{\sqrt{s}}\leq 1-\varepsilon\Big\},\\
\sigma_2&=\Big\{s\in\sigma(A):\,1+\varepsilon\leq \frac{f(s)}{\sqrt{s}}\leq K\Big\},\\
\sigma_3&=\Big\{s\in\sigma(A):\,1-\varepsilon<\frac{f(s)}{\sqrt{s}}<1+\varepsilon\Big\}.
\end{align*}
The proofs of Theorem~\ref{MAIN} and the subsequent Proposition~\ref{propsub}
are based on the following four lemmas,
which hold for all regular initial data, that is,
for all $(u_0,v_0)$ in the domain of $\A$.
Given the trajectory
$$(u(t),\dot u(t))=S(t)(u_0,v_0)\in\D(\A),$$
we split the corresponding energy at time $t\geq 0$ into the sum
$$\E(t)=\sum_{\imath=0}^3\E_\imath(t),
$$
where
$$
\E_\imath (t)= \|A^{\frac12}E_A(\sigma_\imath)u\|^2+
\|E_A(\sigma_\imath)\dot u\|^2.
$$
Note that $\E_\imath \equiv 0$ if $\sigma_\imath=\emptyset$.
The constants $K$ and $\varepsilon$ in the
statements of the lemmas are understood to be the ones
appearing in the definitions of $\sigma_\imath$.

\begin{lemma}
\label{LemmaZero}
For every $K\geq 2$ large enough we have the inequality
$$\E_0 (t)\leq 3\E_0(0)e^{-2m_* t}.
$$
\end{lemma}

\begin{lemma}
\label{LemmaUno}
For every $\varepsilon\in(0,1)$ we have the inequality
$$\E_1 (t)\leq \frac{2-\varepsilon}{\varepsilon}\E_1 (0)e^{-2m_* t}.
$$
\end{lemma}

\begin{lemma}
\label{LemmaDue}
For every $\varepsilon\in(0,1)$ and every $K\geq 2$,
we have the inequality
$$\E_2 (t)\leq\frac{9K^2}\varepsilon \E_2 (0)e^{-2m_* t}.
$$
\end{lemma}

\begin{lemma}
\label{LemmaTre}
For every $\varepsilon\in(0,\frac1{16})$ such that $\sigma_3\neq\emptyset$,
we have the inequality
$$\E_3 (t)\leq \frac{8}{\varepsilon}\E_3 (0)e^{-2m_3(1-4\sqrt{\varepsilon}\,) t},
$$
where $m_3=\inf_{s\in \sigma_3} \phi(s)$.
\end{lemma}

The rather technical proofs of the four lemmas, heavily based on the functional calculus of $A$, will be
carried out in the final Sections~\ref{SecL0}, \ref{SecL1}, \ref{SecL2} and~\ref{SecL3} of the paper.

\subsection{Proof of Theorem \ref{MAIN}}
\label{SubProof}
We consider the two cases separately.

\medskip
\noindent
{\bf I.}
Let us choose $K\geq 2$ sufficiently large in order for Lemma~\ref{LemmaZero} to hold.
We claim that, up to fixing $\varepsilon\in(0,\frac1{16})$ small enough,
$$\E_3 (t)\leq \frac8{\varepsilon}\E_3 (0)e^{-2m_* t}.$$
This is trivially true if
$\sigma_3=\emptyset$ for some $\varepsilon$.
Let us then assume that $\sigma_3\neq\emptyset$ for all $\varepsilon\in (0,\frac1{16})$.
We preliminarily observe that, for all $\varepsilon$ small,
$$m_3>m_*.$$
If not, there would exist a sequence $s_n\in \sigma(A)$ such that
$$\frac{f(s_n)}{\sqrt{s_n}}\to 1\qquad\text{and}\qquad \phi(s_n)\to m_*=m_3.$$
If $s_n$ is unbounded, then (up to a subsequence) $s_n\to\infty$. Therefore,
$$m_*\sim\phi(s_n)\sim f(s_n)\sim\sqrt{s_n}\,\to\infty,$$
a contradiction. This tells that $s_n$ must be bounded. But
then (up to a subsequence) $s_n\to s_*\in\sigma(A)$, as $\sigma(A)$ is closed.
By continuity, we deduce that
$f(s_*)=\sqrt{s_*}$ and $\phi(s_*)=m_*$, against the assumption that $S(t)$ is not resonant. At this point,
the claim follows by fixing $\varepsilon$ in Lemma~\ref{LemmaTre}
small enough that $m_3(1-4\sqrt{\varepsilon}\,)\geq m_*$
(just note that $m_3$ is a decreasing function of $\varepsilon$).
Once the claim is proven, setting
$$M=\frac{9K^2}{\varepsilon},$$
and  taking advantage of Lemmas~\ref{LemmaZero}-\ref{LemmaDue}, we conclude
that
$$\E (t)=\sum_{\imath=0}^3\E_\imath(t)\leq M\E(0)e^{-2m_* t}.
$$
Due to the continuity of the semigroup,
the latter inequality remains valid by density for all initial data $(u_0,v_0)$ belonging
to the phase space $\H$. This finishes the proof.
\qed

\medskip
\noindent
{\bf II.} In this case, we have the equality $m_3=m_*$ for every $\varepsilon\in(0,\frac1{16})$.
Select then any $K\geq 2$ sufficiently large in order for Lemma~\ref{LemmaZero} to hold. From
Lemmas~\ref{LemmaZero}-\ref{LemmaTre} we learn that the inequality
$$\E (t)\leq \frac{9K^2}{\varepsilon}\E(0)e^{-2m_*(1-4\sqrt{\varepsilon}\,)t}
$$
holds for every $\varepsilon\in(0,\frac1{16})$.
Fixing an arbitrary $\varepsilon_*\in (0,\frac1{16})$,
we choose for every $t\geq 0$
$$\varepsilon=\varepsilon(t)=\frac{\varepsilon_*}{(1+t)^2}.$$
Then, calling
$$M=\frac{9K^2}{\varepsilon_*}e^{8m_* \sqrt{\varepsilon_*}\,},
$$
we finally obtain
$$\E (t)\leq M(1+t)^2\E(0)e^{-2m_* t}.
$$
As before, the sought inequality is valid by density for all
initial data.
\qed

\subsection{Proof of Proposition \ref{propsub}}
If the semigroup is subdamped,  then $\sigma(A)=\sigma_1$ with $\varepsilon=1-\ell$.
Accordingly, the energy $\E$ reduces to $\E_1$, and point {\bf I} of Theorem~\ref{MAIN}
is nothing but Lemma~\ref{LemmaUno}.
\qed
%%%%%%%%%%%%%%%%%%%%%%%%%%%%%%%%%%%%%%%%%%%%

%%%%%%%%%%%%%%%%%%%%%%%%%%%%%%%%%%%%%%%%%%%%
\section{An Application:\ Wave Equations with Fractional Damping}
\label{Sec5}

\noindent
For $a>0$ and $\theta\in[0,1]$, we consider equation~\eqref{WAVE} with
$f(s)=a s^\theta$,
namely,
$$
\ddot u + 2 a A^\theta \dot u +A u=0.
$$
The function $f$ clearly complies with~\eqref{STRUCT}, so that Theorem~\ref{MAIN}
and Proposition~\ref{propsub} apply.
We denote by $s_0>0$ the minimum
of the spectrum of $A$.

\begin{remark}
Given a bounded domain $\Omega\subset \R^N$
with smooth boundary $\partial \Omega$, a concrete instance of this model is obtained by taking
$H=L^2(\Omega)$ and the Laplace-Dirichlet operator $A=-\Delta$ with domain
$$\D(-\Delta)=H^2(\Omega)\cap H_0^1(\Omega).
$$
In particular, for $\theta=0$ we have the weakly damped wave equation,
whereas for $\theta=1$ we have the strongly damped wave equation.
Another physically relevant model is the beam (for $N=1$) or plate (for $N=2$) equation with fractional damping,
subject to the hinged boundary conditions.
In this case, the operator $A$ is the Bilaplacian $\Delta^2$ with domain
$$\D(\Delta^2)=\big\{u\in H^2(\Omega)\cap H_0^1(\Omega):\, \Delta u\in
H^2(\Omega)\cap H_0^1(\Omega)\big\}.
$$
\end{remark}

For this choice of $f$, the function $\phi$ reads
$$
\phi(s)=\begin{cases}
as^\theta &\quad\text{if } s\geq a^{\frac{2}{1-2\theta}},\\
\noalign{\vskip1.5mm}
as^\theta -\sqrt{a^2s^{2\theta}-s} &\quad\text{if } s< a^{\frac{2}{1-2\theta}},
\end{cases}
\qquad\text{if } \theta<\frac12,
$$
and
$$
\phi(s)=\begin{cases}
as^\theta &\quad\text{if } s\leq a^{\frac{2}{1-2\theta}},\\
\noalign{\vskip1.5mm}
as^\theta -\sqrt{a^2s^{2\theta}-s} &\quad\text{if } s> a^{\frac{2}{1-2\theta}},
\end{cases}
\qquad\text{if } \theta>\frac12,
$$
whereas
$$
\phi(s)=\begin{cases}
a\sqrt{s}\,&\quad\text{if } a\leq 1,\\
\noalign{\vskip1.5mm}
\big(a-\sqrt{a^2-1}\,\big)\sqrt{s}\,&\quad\text{if } a> 1,
\end{cases}
\qquad\text{if } \theta=\frac12.
$$
Note that $\phi(s)\leq \sqrt{s}\,$. Besides, for $\theta\neq \frac12$, equality occurs
only for $s=a^{\frac{2}{1-2\theta}}$. Accordingly, whenever $\theta\neq\frac12$
and the semigroup is resonant,
$$m_*=a^{\frac{1}{1-2\theta}}.
$$
In general, the picture strongly depends on the exponent $\theta$.
To this end, we first observe that
$$\frac{f(s)}{\sqrt{s}}\,=as^{\theta-\frac12}$$
is decreasing for $\theta<\frac12$, constant and equal to $a$
for $\theta=\frac12$,  increasing and diverging to infinity for $\theta>\frac12$.
In particular, this tells that, if $\theta>\frac12$ and
$A$ is an unbounded operator, the semigroup $S(t)$ is never subdamped.

\begin{figure}
\centering
\begin{subfigure}{.4\linewidth}
\includegraphics[width=\linewidth]{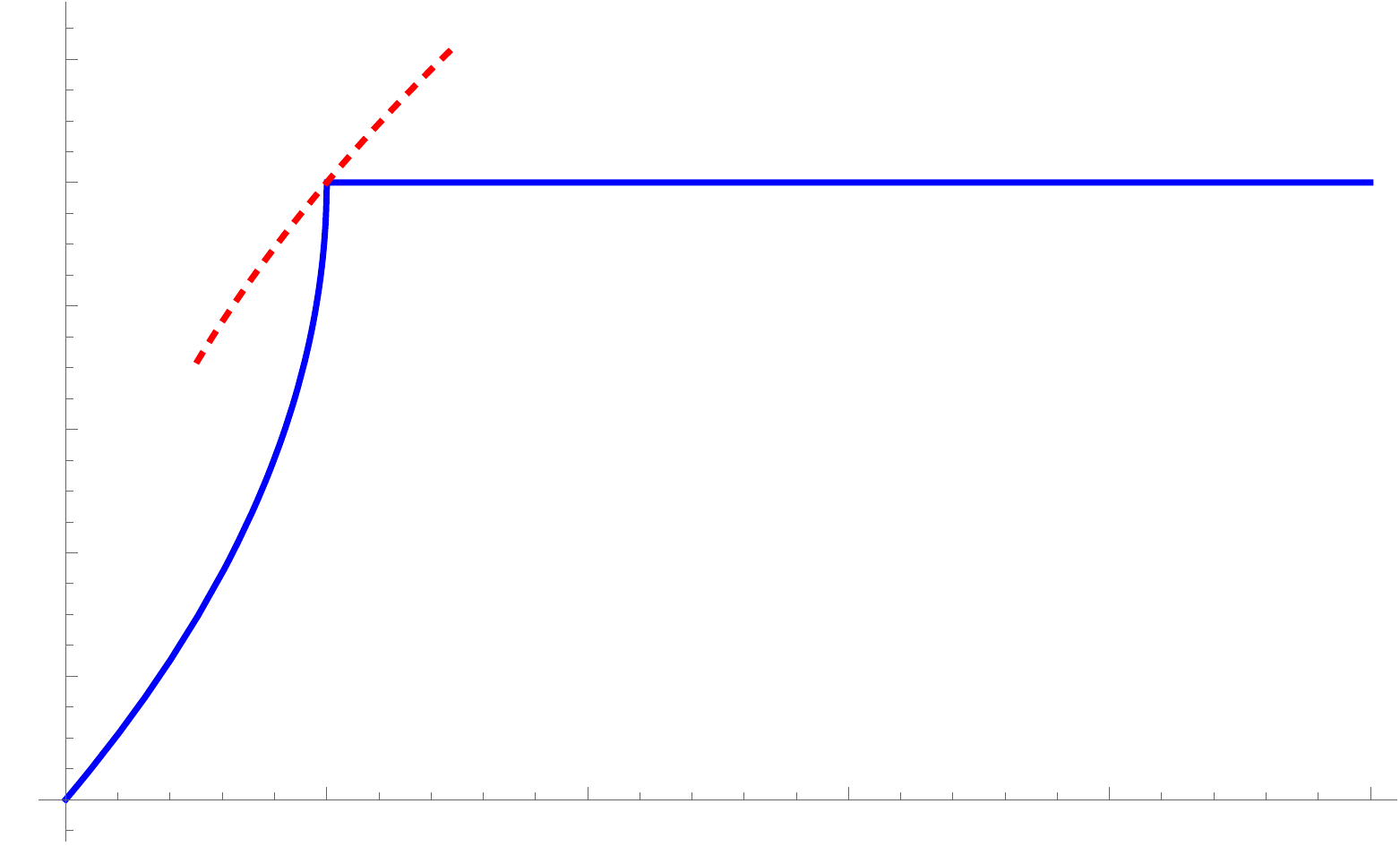}
\caption*{\footnotesize$\boldsymbol{\theta = 0}$}\label{fig1}
\end{subfigure}\qquad
\begin{subfigure}{.4\linewidth}
\includegraphics[width=\linewidth]{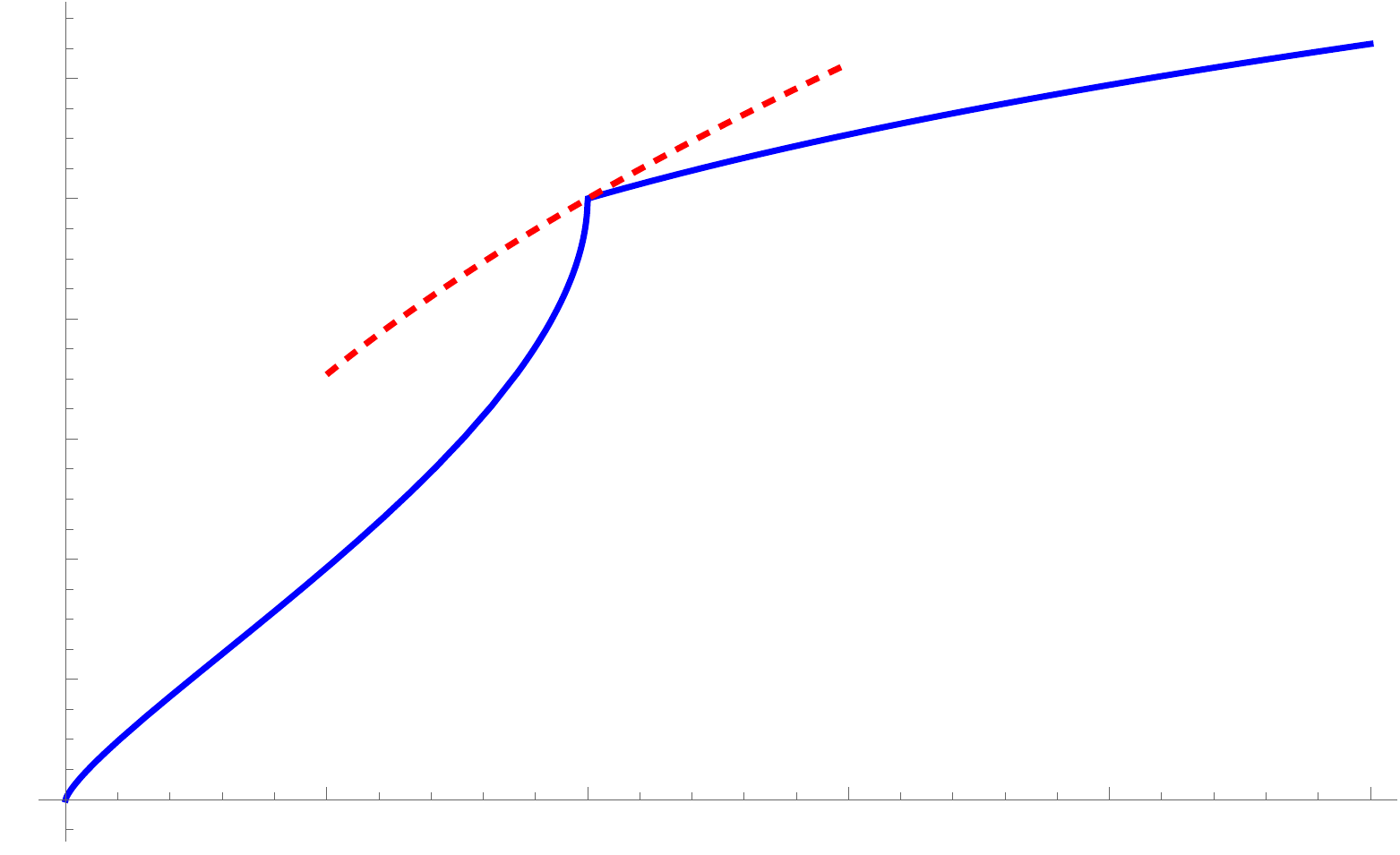}
\caption*{\footnotesize$\boldsymbol{\theta \in(0,\frac12)}$}\label{fig2}
\end{subfigure}
\begin{subfigure}{.4\linewidth}
\includegraphics[width=\linewidth]{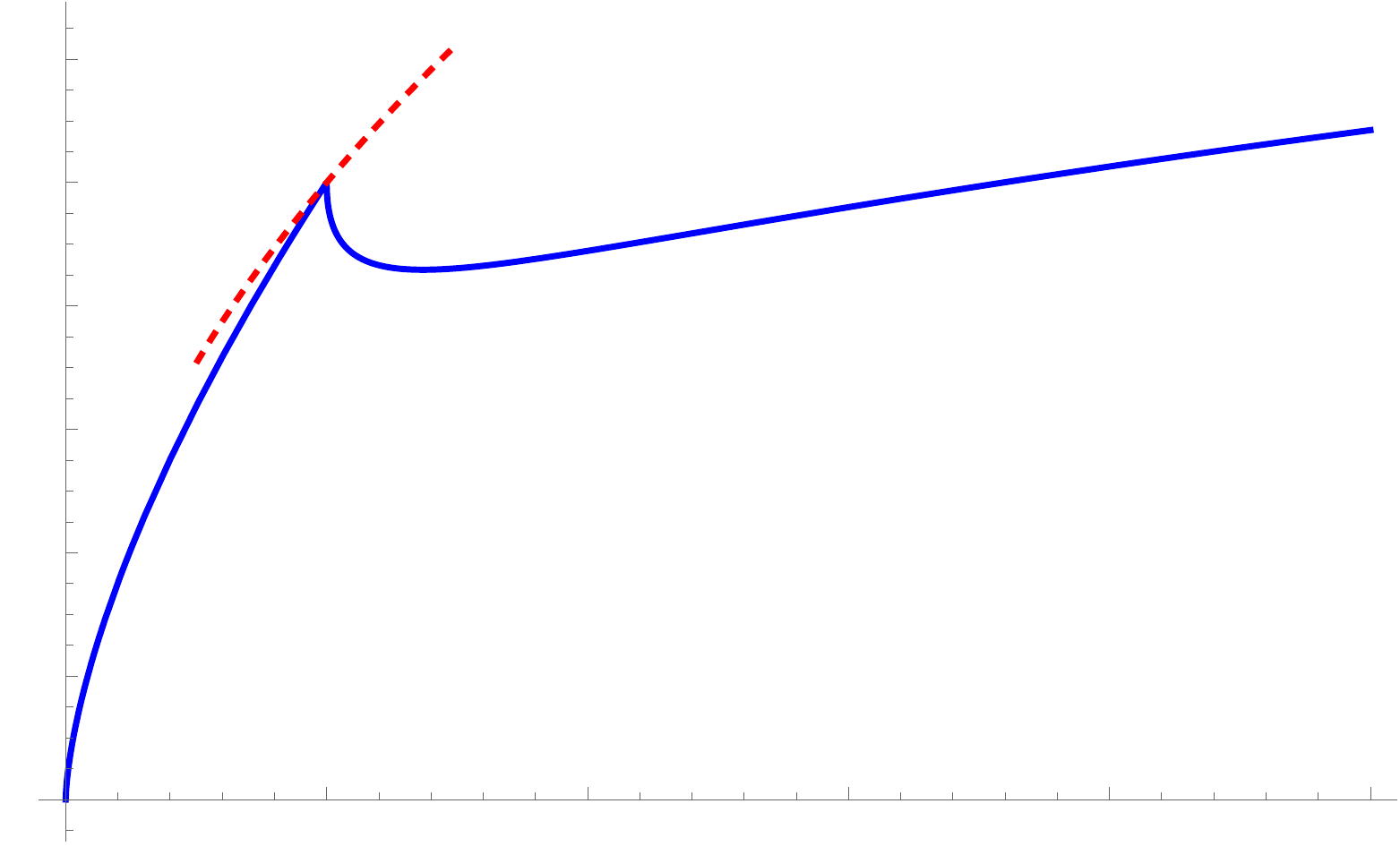}
\caption*{\footnotesize$\boldsymbol{\theta \in(\frac12,1)}$}\label{fig3}
\end{subfigure}\qquad
\begin{subfigure}{.4\linewidth}
\includegraphics[width=\linewidth]{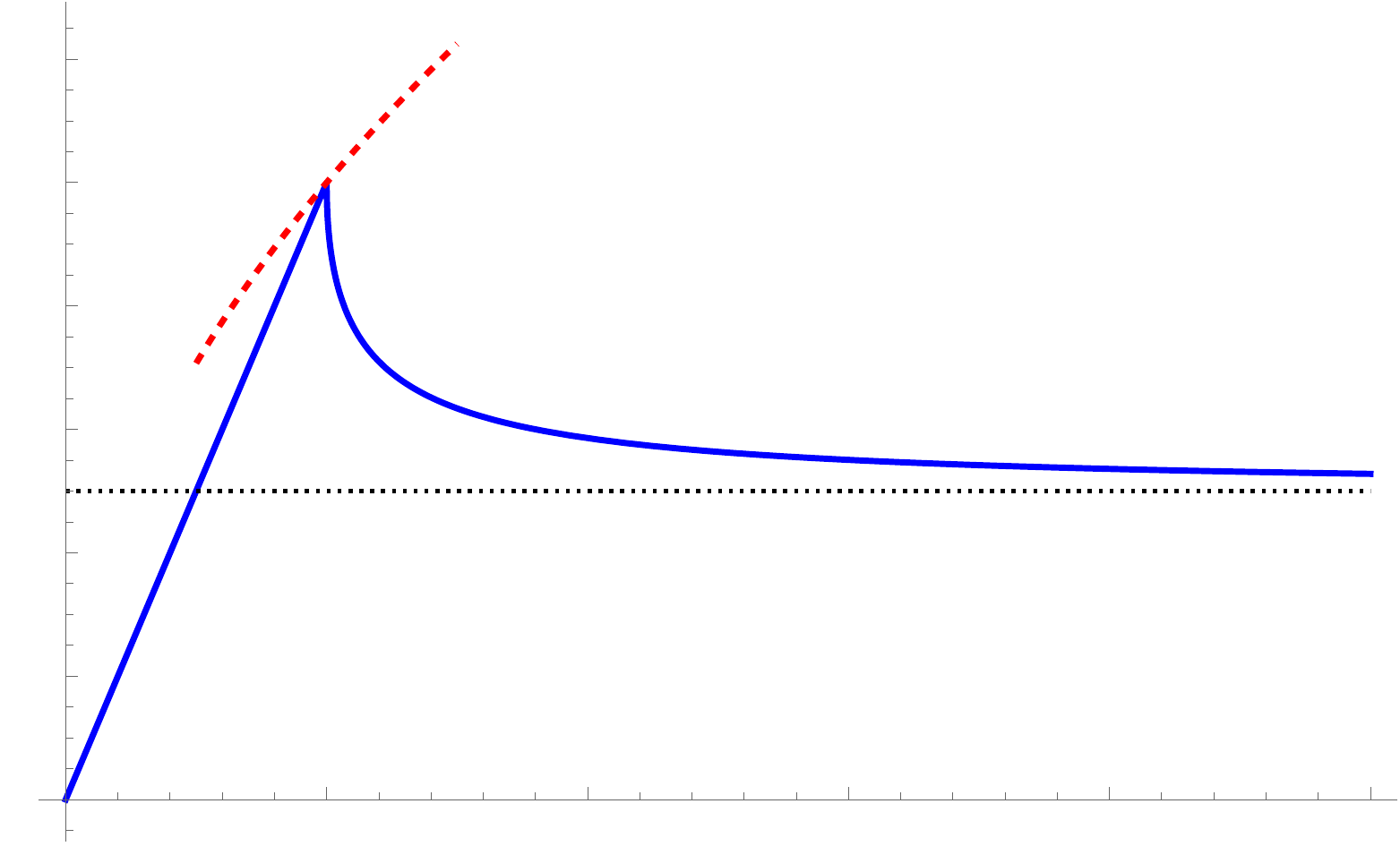}
\caption*{\footnotesize$\boldsymbol{\theta = 1}$}\label{fig4}
\end{subfigure}
\caption{Profile of $\phi(s)$ versus
$\sqrt{s}$ (dashed) for different values of $\theta$.}
\label{fig:all}
\end{figure}

\subsection*{$\bullet$ Case $\boldsymbol{\theta\in[0,\tfrac12)}$}
The function $\phi$ is increasing (actually, strictly increasing if $\theta>0$).
Therefore, Theorem~\ref{MAIN} holds with
$$m_*=\phi(s_0).$$
The semigroup $S(t)$ is resonant if and only if $s_0=a^{\frac{2}{1-2\theta}}$.
Besides, the semigroup is subdamped whenever $s_0>a^{\frac{2}{1-2\theta}}$. In which case,
we infer from Proposition~\ref{propsub} that the constant $C$
of the theorem reads
$$C=\sqrt{\frac{{s_0}^{\frac12-\theta}+a}{{s_0}^{\frac12-\theta}-a}}\,.
$$

\subsection*{$\bullet$ Case $\boldsymbol{\theta=\frac12}$} Again, $m_*=\phi(s_0)$.
The semigroup
is resonant if and only if $a=1$, and is subdamped if and only if $a<1$, with
$$C=\sqrt{\frac{1+a}{1-a}}\,.
$$
Note that $C$ is now independent of $s_0$.

\subsection*{$\bullet$ Case $\boldsymbol{\theta\in(\tfrac12,1)}$}
This is the most intriguing situation, since:
\begin{itemize}
\item[{\bf -}] $\phi$ is increasing for $s<a^{\frac{2}{1-2\theta}}$.
\smallskip
\item[{\bf -}] $\phi$ is decreasing for $s\in(a^{\frac{2}{1-2\theta}}, s_{\rm m})$, where
$s_{\rm m}>a^{\frac{2}{1-2\theta}}$ is the minimum of $\phi$,
here thought as defined on $(0,\infty)$ and not only on $\sigma(A)$,
on the interval $(a^{\frac{2}{1-2\theta}},\infty)$.
\smallskip
\item[{\bf -}] $\phi$ is increasing and diverging to infinity for $s>s_{\rm m}$.
\end{itemize}
The value $m_*$ is now
less explicit, since it might not be equal to $\phi(s_0)$.
In particular, resonance occurs if and only if
$$s_0=a^{\frac{2}{1-2\theta}}\qquad\text{and}\qquad
\sigma(A)\cap (a^{\frac{2}{1-2\theta}},s_{\rm b})=\emptyset,
$$
where $s_{\rm b}$ is the bigger of the two (distinct) solutions to the equation
$$s-2a^{\frac{2-2\theta}{1-2\theta}}s^{\theta}+a^{\frac{2}{1-2\theta}}=0,\quad s>0,
$$
the other one being $a^{\frac{2}{1-2\theta}}$. In this case, $s_0$ is an eigenvalue
of $A$, being an isolated point of $\sigma(A)$, and in turn $-m_*$ is an eigenvalue of $\A$.

\subsection*{$\bullet$ Case $\boldsymbol{\theta=1}$} We have the equality $\phi(s)=as$  for
$s<\frac1{a^2}$, where $\phi$,
thought again as defined on $(0,\infty)$, reaches its maximum value, equal to $\frac1a$. After that, $\phi$ is decreasing,
and converges to $\frac1{2a}$ as $s\to\infty$. Then, if $A$ is unbounded,
$$m_*=\min\Big\{as_0,\frac1{2a}\Big\},$$
whereas if $A$ is bounded,
$$m_*=\min\big\{as_0,\phi(s_{\rm M})\big\},$$
where $s_{\rm M}$ is the maximum of $\sigma(A)$.
Quite interestingly, when $\theta=1$ resonance cannot occur, except in the trivial situation
where $\sigma(A)=\{\frac1{a^2}\}$.
%%%%%%%%%%%%%%%%%%%%%%%%%%%%%%%%%%%%%%%%%%%%

%%%%%%%%%%%%%%%%%%%%%%%%%%%%%%%%%%%%%%%%%%%%
\section{Proof of  Lemma \ref{LemmaZero}}
\label{SecL0}

\noindent
We choose $K$ sufficiently large in order to have the set inclusion
\begin{equation}
\label{ZeroTRE}
\sigma_0\subset (K,\infty).
\end{equation}
Indeed, recalling the second constraint in~\eqref{STRUCT}, if $s\leq K$ with $K$ large then
$$\frac{f(s)}{\sqrt{s}}\leq \sup_{s\in\sigma(A)}\frac{f(s)}{s}\,\sqrt{K}\,\leq K.
$$
Defining the function
$$
g(s) = f(s) - \frac{m_*}{2},
$$
we rewrite equation \eqref{WAVE} in the equivalent form
$$
\ddot u + 2 g(A) \dot u + m_* \dot u + A u = 0.
$$
Note that, due to~\eqref{ossbella},
\begin{equation}
\label{disugG}
s\geq 2m_*g(s),\quad \forall s \in \sigma_0.
\end{equation}
Next, setting $w=E_A(\sigma_0)u$, we introduce the functional
$$\F_0=\E_0 + 2m_*^2\|w\|^2 + 4 m_* \langle w,\dot w\rangle.$$
In what follows, as well as in the proofs of the remaining three lemmas, all the calculations are rigorous,
since we are working with regular solutions. In particular, we use the fact that the selfadjoint projection
$E_A(\sigma_0)$ commutes with all the functions of $A$.
Multiplying the equation by the test function
$2\dot w+4m_*w$, we are led to the differential identity
$$\frac{d}{dt}\F_0+2m_*\F_0+2\G_0=0,
$$
where
$$\G_0=m_*\big[\|A^{\frac12}w\|^2-\|(2m_*g)^{\frac12}(A)w\|^2\big]
+2\big[\|g^{\frac12}(A)(m_*w+\dot w)\|^2-m_*\|m_*w+\dot w\|^2\big].$$
In light of~\eqref{disugG}, we see at once that
$$\|A^{\frac12}w\|^2-\|(2m_*g)^{\frac12}(A)w\|^2\geq 0.$$
Besides, from the very definition of $g$ and $\sigma_0$,
$$
g(s) \geq K \sqrt{s} - \frac{m_*}{2}, \quad \forall s \in \sigma_0.
$$
Thus, we infer from~\eqref{ZeroTRE} that $g(s) \geq m_*$ provided that $K$ is large enough.
Accordingly,
$$\|g^{\frac12}(A)(m_*w+\dot w)\|^2-m_*\|m_*w+\dot w\|^2\geq 0.$$
In summary, we have proved that
$\G_0\geq 0$. Therefore,
$$
\frac{d}{dt}\F_0+2m_*\F_0\leq 0,
$$
and by the Gronwall Lemma we conclude that
$$\F_0(t)\leq \F_0(0)e^{-2m_* t}.$$
The proof is finished once we show the double inequality
$$
\frac12 \E_0 \leq \F_0 \leq \frac32\E_0.
$$
Indeed, from the Cauchy-Schwarz and the Young inequalities,
and exploiting~\eqref{ZeroTRE},
$$
2m_*^2\|w\|^2+4 m_* |\langle w,\dot w\rangle|
\leq 10m_*^2 \|w\|^2+\frac12 \|\dot w\|^2
\leq\frac{10m_*^2}{K} \|A^{\frac12}w\|^2+\frac12 \|\dot w\|^2\leq \frac12\E_0,
$$
as soon as $K$ is sufficiently large.
\qed
%%%%%%%%%%%%%%%%%%%%%%%%%%%%%%%%%%%%%%%%%%%%

%%%%%%%%%%%%%%%%%%%%%%%%%%%%%%%%%%%%%%%%%%%%
\section{Proof of  Lemma \ref{LemmaUno}}
\label{SecL1}

\noindent
Setting $w=E_A(\sigma_1)u$, we introduce the functional
$$\F_1=\E_1
+ 2\langle f(A)w,\dot w\rangle.$$
Multiplying~\eqref{WAVE} by the test function
$2\dot w + 2 f(A)w$,
we get the differential identity
$$\frac{d}{dt}\F_1+2m_*\F_1+2\G_1=0,
$$
having set (recall that in this region $s-f^2(s)>0$)
\begin{align*}
\G_1 &= \|f^{\frac12}(A)(A-f^2(A))^{\frac12} w\|^2 - m_*\|(A-f^2(A))^{\frac12} w\|^2 \\
&\quad+\|f^{\frac12}(A)(f(A)w + \dot w)\|^2 - m_* \|f(A)w + \dot w\|^2.
\end{align*}
Since
$$m_*\leq \phi(s)=f(s),\quad\forall s\in\sigma_1,$$
it is apparent that $\G_1\geq 0$.
Therefore, we arrive at
\begin{equation}
\label{diffineq}
\frac{d}{dt}\F_1+2m_*\F_1\leq 0.
\end{equation}
Finally, we claim that
\begin{equation}
\label{equiv1}
\varepsilon \E_1 \leq \F_1 \leq (2-\varepsilon) \E_1.
\end{equation}
Indeed, recalling that $f(s)\leq (1-\varepsilon)\sqrt{s}$ for every $s\in \sigma_1$, from the Cauchy-Schwarz and the Young inequalities,
we get
$$
2|\langle f(A)w,\dot w\rangle|
\leq (1-\varepsilon)\|A^{\frac12}w\|^2
+(1-\varepsilon)\|\dot w\|^2.
$$
Collecting \eqref{diffineq}-\eqref{equiv1}, and exploiting the Gronwall
Lemma, the conclusion follows.
\qed
%%%%%%%%%%%%%%%%%%%%%%%%%%%%%%%%%%%%%%%%%%%%

%%%%%%%%%%%%%%%%%%%%%%%%%%%%%%%%%%%%%%%%%%%%
\section{Proof of  Lemma \ref{LemmaDue}}
\label{SecL2}

\noindent
In this proof, as well as in the following one, the symbol $\langle\cdot,\cdot\rangle$ will also be
used to denote the duality product.
Setting $w=E_A(\sigma_2)u$, we introduce the functional
$$\F_2=\E_2
+2\langle (f^2(A)-A)w,w\rangle+ 2\langle f(A)w,\dot w\rangle.$$
Multiplying~\eqref{WAVE} by the test function $2\dot w+2f(A)w$
we get the differential identity
$$\frac{d}{dt}\F_2+2m_*\F_2+2\G_2=0,
$$
having set
\begin{align*}
\G_2& = \langle (Af(A)+m_*A-2m_*f^2(A)) w,w\rangle+
\langle (f(A)-m_*)\dot w,\dot w\rangle\\
&\quad+2\langle (A-m_*f(A)) w,\dot w\rangle.
\end{align*}
We claim that $\G_2\geq 0$. Indeed, since
$$m_*\leq \phi(s)=f(s)-\sqrt{f^2(s)-s}\,< f(s),\quad\forall s\in\sigma_2,$$
exploiting the Cauchy-Schwarz and the Young inequalities,
\begin{align*}
&2\langle (A-m_*f(A)) w,\dot w\rangle\\
&=2\langle (f(A)-m_*)^{-\frac12}(A-m_*f(A)) w,(f(A)-m_*)^{\frac12}\dot w\rangle\\
&\geq -\langle (f(A)-m_*)^{-1}(A-m_*f(A))^2w, w\rangle
-\langle (f(A)-m_*)\dot w,\dot w\rangle.
\end{align*}
Therefore,
since
$$sf(s)+m_*s-2m_*f^2(s)-\frac{(s-m_*f(s))^2}{f(s)-m_*}=\frac{f^2(s)-s}{f(s)-m_*}(s+m_*^2-2m_*f(s)),
$$
we get
$$\G_2\geq \langle (f(A)-m_*)^{-1}(f^2(A)-A)(A+m_*^2-2m_*f(A))w,w\rangle.
$$
The claim follows by noting that, for every $s\in\sigma_2$,
$$s+m_*^2-2m_*f(s)
=\big(f(s)-\sqrt{f^2(s)-s}\,-m_*\big)\big(f(s)+\sqrt{f^2(s)-s}\,-m_*\big)\geq 0.$$
Accordingly, we end up with the differential inequality
\begin{equation}
\label{diffineq2}
\frac{d}{dt}\F_2+2m_*\F_2\leq 0.
\end{equation}
Next, we show that
\begin{equation}
\label{equiv2}
\frac\varepsilon3\E_2 \leq \F_2 \leq 3K^2 \E_2.
\end{equation}
To this end, we further apply
the Cauchy-Schwarz and the Young inequalities, to get
$$2|\langle f(A)w,\dot w\rangle|\leq \frac{1}{\tau}\|f(A)w\|^2
+\tau \|\dot w\|^2,$$
for every $\tau>0$.
Since $f^2(s)\leq K^2s$ for every $s\in\sigma_2$, by setting $\tau=1$ we obtain
$$\F_2\leq \E_2+2(K^2-1)\|A^{\frac12}w\|^2
+K^2\|A^{\frac12}w\|^2
+\|\dot w\|^2\leq 3K^2\E_2.
$$
As far as the other inequality is concerned, we have
$$\F_2\geq (1-\tau)\E_2+\langle((2-\tau^{-1})f^2(A)+(\tau-2)A)w,w\rangle.
$$
Recalling that $f^2(s)\geq (1+\varepsilon)^2s$ for every $s\in \sigma_2$,
and selecting $\tau=1-\frac\varepsilon3$, we end up with
$$\F_2\geq \frac\varepsilon3\E_2+\varrho(\varepsilon)\|A^{\frac12}w\|^2,
$$
where, for every $\varepsilon\in(0,1)$,
$$\varrho(\varepsilon)=\frac{(3-2\varepsilon)(1+\varepsilon)^2}{3-\varepsilon}-1-\frac{\varepsilon}3
=\frac{2\varepsilon(6-\varepsilon-3\varepsilon^2)}{3(3-\varepsilon)}>0.
$$
At this point, the conclusion is drawn from
\eqref{diffineq2}-\eqref{equiv2} and the Gronwall
Lemma.
\qed
%%%%%%%%%%%%%%%%%%%%%%%%%%%%%%%%%%%%%%%%%%%%

%%%%%%%%%%%%%%%%%%%%%%%%%%%%%%%%%%%%%%%%%%%%
\section{Proof of  Lemma \ref{LemmaTre}}
\label{SecL3}

\noindent
Along the proof, besides the constraint $\varepsilon\in(0,\frac1{16})$,
we will use the fact that
$$m_3\leq \phi(s)\leq f(s)< (1+\varepsilon) \sqrt{s}\,,\quad\forall s\in\sigma_3.$$
Setting $w=E_A(\sigma_3)u$, we introduce the functional
$$
\F_3 = \E_3+ 2 (1+4\sqrt{\varepsilon}\,) \langle f(A) w,\dot w\rangle
+ 8\sqrt{\varepsilon}\, (1+4\sqrt{\varepsilon}\,) \|f(A)w\|^2.
$$
Multiplying~\eqref{WAVE} by the test function $2\dot w+ 2(1+4\sqrt{\varepsilon}\,)f(A)w$, we get the
differential identity
$$
\frac{d}{dt} \F_3+ 2m_3(1-4\sqrt{\varepsilon}\,)\F_3 + 2 \G_3=0.
$$
Here,
\begin{align*}
\G_3 = \langle p(A)w,w\rangle+(1-4\sqrt{\varepsilon}\,)\langle (f(A)-m_3)\dot w,\dot w\rangle
+2(1-16\varepsilon)\langle f(A)(f(A)-m_3)w,\dot w\rangle,
\end{align*}
having set
\begin{align*}
p(s)&=(1+4\sqrt{\varepsilon}\,)s f(s) - m_3(1-4\sqrt{\varepsilon}\,) s - 8 m_3 \sqrt{\varepsilon}\, (1-16\varepsilon)f^2(s)\\
&\geq(1+4\sqrt{\varepsilon}\,)s f(s) - m_3(1-4\sqrt{\varepsilon}\,) s - 8 m_3 \sqrt{\varepsilon}\, (1-16\varepsilon)(1+\varepsilon)^2 s,\\
&=(1+4\sqrt{\varepsilon}\,)s (f(s)-m_3)
+8 m_3\sqrt{\varepsilon}\,\big[1-(1-16\varepsilon)(1+\varepsilon)^2\big]s\\
&\geq(1+4\sqrt{\varepsilon}\,)s (f(s)-m_3).
\end{align*}
Accordingly,
\begin{align*}
\G_3 &\geq  (1+4\sqrt{\varepsilon}\,) \langle A(f(A)-m_3)w,w\rangle+(1-4\sqrt{\varepsilon}\,)\langle (f(A)-m_3)\dot w,\dot w\rangle\\
&\quad +2(1-16\varepsilon)\langle f(A)(f(A)-m_3)w,\dot w\rangle.
\end{align*}
By an application of the
Cauchy-Schwarz and the Young inequalities, we find
\begin{align*}
&2(1-16\varepsilon)\langle f(A)(f(A)-m_3)w,\dot w\rangle\\
&=2\langle (1+4\sqrt{\varepsilon}\,)(1-4\sqrt{\varepsilon}\,)^{\frac12} f(A)(f(A)-m_3)^{\frac12}w,(1-4\sqrt{\varepsilon}\,)^{\frac12} (f(A)-m_3)^{\frac12}\dot w\rangle\\
&\geq -(1-16\varepsilon)(1+4\sqrt{\varepsilon}\,)\langle f^2(A)(f(A)-m_3)w ,w\rangle
-(1-4\sqrt{\varepsilon}\,)\langle (f(A)-m_3)\dot w,\dot w\rangle,
\end{align*}
implying in turn
\begin{align*}
\G_3 &\geq  (1+4\sqrt{\varepsilon}\,) \langle (f(A)-m_3)
(A-(1-16\varepsilon)f^2(A))w,w\rangle.
\end{align*}
Observing that
$$s-(1-16\varepsilon)f^2(s)\geq \big[1-(1-16\varepsilon)(1+\varepsilon)^2\big]s\geq 0,
$$
we conclude that
$\G_3\geq 0$, and we arrive at
$$
\frac{d}{dt} \F_3+ 2m_3(1-4\sqrt{\varepsilon}\,)\F_3 \leq 0.
$$
After an application of the Gronwall Lemma, the proof is finished
once we show that
$$
\varepsilon\E_3 \leq \F_3 \leq 8 \E_3.
$$
Indeed, the second inequality is a straightforward consequence
of the Cauchy-Schwarz and the Young inequalities, and is left to the reader.
Concerning the first one,
we have
$$2 (1+4\sqrt{\varepsilon}\,) \langle f(A) w,\dot w\rangle\geq -(1-\varepsilon)\|\dot w\|^2
-\frac{(1+4\sqrt{\varepsilon}\,)^2}{1-\varepsilon}\|f(A)w\|^2.
$$
Therefore,
\begin{align*}
\F_3&\geq \varepsilon\E_3 +(1-\varepsilon)\|A^{\frac12}w\|^2
+\Big[ 8\sqrt{\varepsilon}\, (1+4\sqrt{\varepsilon}\,)
-\frac{(1+4\sqrt{\varepsilon}\,)^2}{1-\varepsilon}\Big]\|f(A)w\|^2\\
&\geq \varepsilon\E_3
+\Big[\frac{1-\varepsilon}{(1+\varepsilon)^2}+ 8\sqrt{\varepsilon}\, (1+4\sqrt{\varepsilon}\,)
-\frac{(1+4\sqrt{\varepsilon}\,)^2}{1-\varepsilon}\Big]\|f(A)w\|^2,
\end{align*}
and the quantity in the square brackets is positive for all $\varepsilon\in(0,\frac1{16})$.
\qed
%%%%%%%%%%%%%%%%%%%%%%%%%%%%%%%%%%%%%%%%%%%%

%%%%%%%%%%%%%%%%%%%%%%%%%%%%%%%%%%%%%%%%%%%%

%%%%%%%%%%%%%%%%%%%%%%%%%%%%%%%%%%%%%%%%%%%%%%%%%%


\begin{thebibliography}{99}

\bibitem{BATK}
{\au A. B\'atkai, K.-J. Engel},
{\ti Exponential decay of $2\times2$ operator matrix semigroups},
{\jou J.\ Comput.\ Anal.\ Appl.}
\no{6}{153--163}{2004}

\bibitem{CR}
{\au G. Chen, D.L. Russell},
{\it A mathematical model for linear elastic systems with structural damping},
{\jou Quart.\ Appl.\ Math.}
\no{39}{433--454}{1981/82}

\bibitem{CT1}
{\au S. Chen, R. Triggiani},
{\it Proof of extensions of two conjectures on structural damping for elastic systems},
{\jou Pacific J.\ Math.}
\no{136}{15--55}{1989}

\bibitem{CT2}
{\au S. Chen, R. Triggiani},
{\it Gevrey class semigroups arising from elastic systems with gentle
dissipation:\ the case $0<\alpha<\tfrac12$},
{\jou Proc.\ Amer.\ Math.\ Soc.}
\no{110}{401--415}{1990}

\bibitem{CT3}
{\au S. Chen, R. Triggiani},
{\it Characterization of domains of fractional powers of certain
operators arising in elastic systems, and applications},
{\jou J.\ Differential Equations}
\no{88}{279--293}{1990}

\bibitem{SILVER}
{\au M. Conti, L. Liverani, V. Pata},
{\it On the optimal decay rate of the weakly damped wave equation},
{\jou Commun.\ Pure Appl.\ Anal.}
\no{21}{3421--3424}{2022}

\bibitem{VDRM}
{\au M. Conti, F. Dell'Oro, V. Pata},
{\ti Some unexplored questions arising in linear viscoelasticity},
{\jou J.\ Funct.\ Anal.}
\no{282}{Paper No.\ 109422, 43 pp}{2022}

\bibitem{DPWAVE}
{\au F. Dell'Oro, V. Pata},
{\ti Second order linear evolution equations with general dissipation},
{\jou Appl.\ Math.\ Optim.}
\no{83}{1877--1917}{2021}

\bibitem{ENG}
{\au K.-J. Engel, R. Nagel},
{\bk One-parameter semigroups for linear evolution equations},
\eds{Springer-Verlag}{New York}{2000}

\bibitem{GOLD1}
{\au G.R. Goldstein, J.A. Goldstein, G. Perla Menzala},
{\ti On the overdamping phenomenon:\ a general result and applications},
{\jou Quart.\ Appl.\ Math.}
\no{71}{183--199}{2013}

\bibitem{GOLD2}
{\au G.R. Goldstein, J.A. Goldstein, G. Reyes},
{\ti Overdamping and energy decay for abstract wave equations with strong damping},
{\jou Asymptot.\ Anal.}
\no{88}{217--232}{2014}

\bibitem{GRIN1}
{\au R.O. Griniv, A.A. Shkalikov},
{\ti Exponential stability of semigroups associated with some operator models in mechanics. (Russian)},
{\jou translation in Math.\ Notes}
\no{73}{618--624}{2003}

\bibitem{HUA1}
{\au F. Huang},
{\ti On the holomorphic property of the semigroup associated
with linear elastic systems with structural damping},
{\jou Acta Math.\ Sci.\ (English Ed.)}
\no{5}{271--277}{1985}

\bibitem{HUA2}
{\au F. Huang},
{\ti On the mathematical model for linear elastic systems with analytic damping},
{\jou SIAM J.\ Control Optim.}
\no{26}{714--724}{1988}

\bibitem{HUA3}
{\au F. Huang, K. Liu},
{\ti Holomorphic property and exponential stability of the semigroup
associated with linear elastic systems with damping},
{\jou Ann.\ Differential Equations}
\no{4}{411--424}{1988}

\bibitem{JT}
{\au B. Jacob, C. Trunk},
{\ti Spectrum and analyticity of semigroups arising in elasticity theory and hydromechanics},
{\jou Semigroup Forum}
\no{79}{79--100}{2009}

\bibitem{LTBOOK}
{\au I. Lasiecka, R. Triggiani},
{\bk Control theory for partial differential equations:\ continuous
and approximation theories},
\eds{Cambridge University Press}{Cambridge}{2000}

\bibitem{LIU}
{\au Z. Liu, J. Yong},
{\ti Qualitative properties of certain $C_0$ semigroups arising in elastic systems with various dampings},
{\jou  Adv. Differential Equations}
\no{3}{643--686}{1998}

\bibitem{PAZ}
{\au A. Pazy},
{\bk Semigroups of linear operators and applications to partial differential equations},
\eds{Springer-Verlag}{New York}{1983}

\bibitem{RUD}
{\au W. Rudin},
{\bk Functional analysis},
\eds{McGraw-Hill}{New York-D\"usseldorf-Johannesburg}{1973}

\end{thebibliography}
\end{document}